\documentclass[a4paper,12pt]{amsart}
\usepackage{amsfonts}
\usepackage{amssymb}
\usepackage{ifthen}
\usepackage{graphicx}
\nonstopmode \numberwithin{equation}{section}

\usepackage[usenames]{color}
\newtheorem{thm}{Theorem}[section]
\newtheorem{lem}{Lemma}[section]
\newtheorem{cor}[thm]{Corollary}
\newtheorem{prop}[thm]{Proposition}

\newtheorem{step}{Step}[section]

\theoremstyle{definition}
\newtheorem{mlem}{Main lemma}[section]
\newtheorem{assertion}{Assertion}[section]
\newtheorem{cl}{Claim}[section]
\newtheorem{ca}{Case}[section]
\newtheorem{sca}{Subcase}[section]
\newtheorem{scl}{Subclaim}[section]
\newtheorem{conj}[thm]{Conjecture}
\newtheorem{fact}{Fact}[section]
\newtheorem{defn}[thm]{Definition}
\newtheorem{op}[thm]{Open Problem}
\newtheorem{prob}{Problem}[section]
\newtheorem{rem}[thm]{Remark}
\newtheorem{exam}[thm]{Example}

\numberwithin{equation}{section}

\newcounter {own}
\def\theown {\thesection       .\arabic{own}}

\newenvironment{pf}[1][]{%
 \vskip 3mm
 \noindent
 \ifthenelse{\equal{#1}{}}%
  {{\slshape Proof. }}%
  {{\slshape #1.} }%
 }%
{\qed\bigskip}

\newcounter{alphabet}
\newcounter{tmp}
\newenvironment{Thm}[1][]{\refstepcounter{alphabet}%
\bigskip%
\noindent%
{\bf Theorem \Alph{alphabet}}%
\ifthenelse{\equal{#1}{}}{}{ (#1)}%
{\bf .} \itshape}{\vskip 8pt}

\makeatletter
\newcommand{\Ref}[1]{\@ifundefined{r@#1}{}{\setcounter{tmp}{\ref{#1}}\Alph{tmp}}}
\makeatother

\newenvironment{Lem}[1][]{\refstepcounter{alphabet}%
\bigskip%
\noindent%
{\bf Lemma \Alph{alphabet}}%
{\bf .} \itshape}{\vskip 8pt}

\newcounter{alphabet2}

\def\be{\begin{equation}}
\def\ee{\end{equation}}

\newcommand{\ben}{\begin{enumerate}}
\newcommand{\een}{\end{enumerate}}

\newcommand{\blem}{\begin{lem}}
\newcommand{\elem}{\end{lem}}
\newcommand{\bthm}{\begin{thm}}
\newcommand{\ethm}{\end{thm}}
\newcommand{\bcor}{\begin{cor}}
\newcommand{\ecor}{\end{cor}}
\newcommand{\beg}{\begin{exam}}
\newcommand{\eeg}{\end{exam}}
\newcommand{\begs}{\begin{examples}}
\newcommand{\eegs}{\end{examples}}
\newcommand{\bdefe}{\begin{defn}}
\newcommand{\edefe}{\end{defn}}
\newcommand{\bques}{\begin{ques}}
\newcommand{\eques}{\end{ques}}
\newcommand{\bei}{\begin{itemize}}
\newcommand{\eei}{\end{itemize}}
\newcommand{\bcon}{\begin{conj}}
\newcommand{\econ}{\end{conj}}
\newcommand{\bop}{\begin{op}}
\newcommand{\eop}{\end{op}}

\newcommand{\bas}{\begin{assertion}}
\newcommand{\eas}{\end{assertion}}

\newcommand{\bfa}{\begin{fact}}
\newcommand{\efa}{\end{fact}}

\newcommand{\bca}{\begin{ca}}
\newcommand{\eca}{\end{ca}}

\newcommand{\bst}{\begin{step}}
\newcommand{\est}{\end{step}}

\newcommand{\bsca}{\begin{sca}}
\newcommand{\esca}{\end{sca}}

\newcommand{\bcl}{\begin{cl}}
\newcommand{\ecl}{\end{cl}}

\newcommand{\bmlem}{\begin{mlem}}
\newcommand{\emlem}{\end{mlem}}

\newcommand{\bscl}{\begin{scl}}
\newcommand{\escl}{\end{scl}}

\newcommand{\bcons}{\begin{conjs}}
\newcommand{\econs}{\end{conjs}}

\newcommand{\bprop}{\begin{prop}}
\newcommand{\eprop}{\end{prop}}

\newcommand{\br}{\begin{rem}}
\newcommand{\er}{\end{rem}}
\newcommand{\brs}{\begin{rems}}
\newcommand{\ers}{\end{rems}}
\newcommand{\bo}{\begin{obser}}
\newcommand{\eo}{\end{obser}}
\newcommand{\bos}{\begin{obsers}}
\newcommand{\eos}{\end{obsers}}
\newcommand{\bpf}{\begin{pf}}
\newcommand{\epf}{\end{pf}}
\newcommand{\ba}{\begin{array}}
\newcommand{\ea}{\end{array}}
\newcommand{\beq}{\begin{eqnarray}}
\newcommand{\beqq}{\begin{eqnarray*}}
\newcommand{\eeq}{\end{eqnarray}}
\newcommand{\eeqq}{\end{eqnarray*}}

\newcommand{\ds}{\displaystyle}

\newcounter{minutes}
\divide\time by 60
\newcounter{hours}
\multiply\time by 60 \addtocounter{minutes}{-\time}

\begin{document}

\bibliographystyle{amsplain}
\title [Harmonic quasiregular mappings]
{Norm estimates of the partial derivatives for harmonic and harmonic elliptic mappings}

\def\thefootnote{}
\footnotetext{ \texttt{\tiny File:~\jobname .tex,
          printed: \number\day-\number\month-\number\year,
          \thehours.\ifnum\theminutes<10{0}\fi\theminutes}
} \makeatletter\def\thefootnote{\@arabic\c@footnote}\makeatother

\author{Shaolin Chen}
 \address{Sh. Chen, College of Mathematics and
Statistics, Hengyang Normal University, Hengyang, Hunan 421008,
People's Republic of China.} \email{mathechen@126.com}

\author{Saminathan Ponnusamy 
}
\address{S. Ponnusamy, Department of Mathematics,
Indian Institute of Technology Madras, Chennai-600 036, India. }
\email{samy@iitm.ac.in}

\author{Xiantao Wang}  
\address{X. Wang, Key Laboratory of High Performance Computing and Stochastic Information Processing,
College of Mathematics and Statistics, Hunan Normal University, Changsha, Hunan 410081, People's Republic of China}
 \email{xtwang@hunnu.edu.cn}



\begin{abstract}
Let $f = P[F]$ denote the Poisson integral of $F$ in the unit disk $\mathbb{D}$ with $F$ being
absolutely continuous in the unit circle $\mathbb{T}$ and $\dot{F}\in L_p(0, 2\pi)$, where $\dot{F}(e^{it})=\frac{d}{dt} F(e^{it})$ and $p\geq 1$.
Recently, the author in \cite{Zhu} proved that $(1)$ if $f$ is a harmonic mapping and $1\leq p< 2$, then $f_{z}$ and $\overline{f_{\overline{z}}}\in \mathcal{B}^{p}(\mathbb{D}),$ the classical Bergman spaces of $\mathbb{D}$ \cite[Theorem 1.2]{Zhu}; $(2)$ if $f$ is a harmonic quasiregular mapping and $1\leq p\leq \infty$, then $f_{z},$ $\overline{f_{\overline{z}}}\in \mathcal{H}^{p}(\mathbb{D}),$ the classical Hardy spaces of $\mathbb{D}$ \cite[Theorem 1.3]{Zhu}. These are the main results in \cite{Zhu}. The purpose of this paper is to generalize these two results.
First, we prove that, under the same assumptions, \cite[Theorem 1.2]{Zhu} is true when $1\leq p< \infty$. Also, we show that \cite[Theorem 1.2]{Zhu} is not true when $p=\infty$. Second, we demonstrate that \cite[Theorem 1.3]{Zhu} still holds true when the assumption $f$ being a harmonic quasiregular mapping is replaced by the weaker one $f$ being a harmonic elliptic mapping.
\end{abstract}

\maketitle \pagestyle{myheadings}
\markboth{ Sh. Chen, S. Ponnusamy and X. Wang}{Harmonic and harmonic elliptic mappings}

\section{Preliminaries and the statement of main results}\label{csw-sec1}

For $a\in\mathbb{C}$ and $r>0$, let $\mathbb{D}(a,r)=\{z:~|z-a|<r\}$. In particular, we use $\mathbb{D}_{r}$ to
denote the disk $\mathbb{D}(0,r)$ and $\mathbb{D}$ to denote the unit disk $\mathbb{D}_{1}$. Moreover, let $\mathbb{T}:=\partial\mathbb{D}$ be the
unit circle.
For $z=x+iy\in\mathbb{C}$, the two complex
differential operators are defined by

$$\frac{\partial}{\partial z}=\frac{1}{2}\left(\frac{\partial}{\partial x}-i\frac{\partial}{\partial y}\right)~\mbox{and}~
\frac{\partial}{\partial \overline{z}}=\frac{1}{2}\left(\frac{\partial}{\partial x}+i\frac{\partial}{\partial y}\right).$$

For $\alpha\in[0,2\pi]$, the {\it directional derivative} of a harmonic
mapping (i.e., a complex-valued harmonic function) $f$ at $z\in\mathbb{D}$ is defined by

$$\partial_{\alpha}f(z)=\lim_{\rho\rightarrow0^{+}}\frac{f(z+\rho e^{i\alpha})-f(z)}{\rho}=f_{z}(z)e^{i\alpha}
+f_{\overline{z}}(z)e^{-i\alpha},$$
where $z+\rho e^{i\alpha}\in\mathbb{D}$, $f_{z}:=\partial f/\partial z$ and $f_{\overline{z}}:=\partial f/\partial \overline{z}$.
Then

$$\|D_{f}(z)\|:=\max\{|\partial_{\alpha}f(z)|:\; \alpha\in[0,2\pi]\}=|f_{z}(z)|+|f_{\overline{z}}(z)|
$$
and
$$l(D_{f}(z)):=\min\{|\partial_{\alpha}f(z)|:\; \alpha\in[0,2\pi]\}=\big||f_{z}(z)|-|f_{\overline{z}}(z)|\big|.
$$

For a sense-preserving harmonic mapping $f$ defined in $\mathbb{D}$,
the Jacobian of $f$ is given by
$$J_{f}=\|D_{f}\|l(D_{f})=|f_{z}|^{2}-|f_{\overline{z}}|^{2},$$ and {\it the second complex dilatation} of $f$ is given by
$\omega=\overline{f_{\overline{z}}}/f_{z}$.  It is well-known that every harmonic
mapping $f$ defined in a simply connected domain $\Omega$ admits a decomposition $f = h + \overline{g}$, where $h$ and $g$ are analytic.
Recall that $f$ is
sense-preserving in $\Omega$ if $J_{f}>0 $ in $\Omega$.
Thus $f$ is locally univalent and sense-preserving in $\Omega$
if and only if $J_{f}>0$ in $\Omega$, which means that $h'\neq
0$ in $\Omega$ and  the analytic function $\omega =g'/h'$ has the
property that  $|\omega (z)|<1$ on $\Omega$ (cf. \cite{Clunie-Small-84,Lewy}).

\subsection*{Hardy type spaces}

For $p\in(0,\infty]$, the {\it generalized Hardy space
$\mathcal{H}^{p}_{\mathcal{G}}(\mathbb{D})$} consists of all
measurable  functions from $\mathbb{D}$ to $\mathbb{C}$ such that
$M_{p}(r,f)$ exists for all $r\in(0,1)$, and $ \|f\|_{p}<\infty$,
where
$$M_{p}(r,f)=\left(\frac{1}{2\pi}\int_{0}^{2\pi}|f(re^{i\theta})|^{p}\,d\theta\right)^{\frac{1}{p}}
$$
and
$$\|f\|_{p}=
\begin{cases}
\displaystyle\sup \{M_{p}(r,f):\; 0<r <1\}
& \mbox{if } p\in(0,\infty),\\
\displaystyle\sup\{|f(z)|:\; z\in\mathbb{D}\} &\mbox{if } p=\infty.
\end{cases}$$

The classical {\it Hardy space $\mathcal{H}^{p}(\mathbb{D})$}, that is, all the elements are analytic, is a subspace of $\mathcal{H}^{p}_{\mathcal{G}}(\mathbb{D})$ (cf. \cite{CPR,Du1}).

\subsection*{Bergman type spaces}
For $p\in(0,\infty]$, the {\it generalized Bergman space
$\mathcal{B}^{p}_{\mathcal{G}}(\mathbb{D})$} consists of all measurable functions $f:\;\mathbb{D}\rightarrow\mathbb{C}$ such that
$$\|f\|_{b^{p}}=
\begin{cases}
\displaystyle\left(\int_{\mathbb{D}}|f(z)|^{p}d\sigma(z)\right)^{\frac{1}{p}}
& \mbox{if } p\in(0,\infty),\\
\displaystyle \mbox{ess}\sup\{|f(z)|:\; z\in \mathbb{D}\} &\mbox{if } p=\infty,
\end{cases}
$$
where  $d\sigma(z)=\frac{1}{\pi}dxdy$ denotes  the normalized Lebesgue area measure on $\mathbb{D}$.
The classical {\it Bergman space $\mathcal{B}^{p}(\mathbb{D})$}, that is, all the elements are analytic, is a subspace of $\mathcal{B}^{p}_{\mathcal{G}}(\mathbb{D})$ (cf. \cite{HKZ}).
Obviously, $\mathcal{H}^{p}(\mathbb{D})\subset\mathcal{B}^{p}(\mathbb{D})$ for each $p\in(0,\infty]$.

\subsection*{Poisson integrals}

Denote by $L^{p}(\mathbb{T})~(p\in[1,\infty])$ the space of all measurable functions  $F$
of $\mathbb{T}$ into $\mathbb{C}$ with

$$\|F\|_{L^{p}}=
\begin{cases}
\displaystyle\left(\frac{1}{2\pi}\int_{0}^{2\pi}|F(e^{i\theta})|^{p}d\theta\right)^{\frac{1}{p}}
& \mbox{if } p\in[1,\infty),\\
\displaystyle \mbox{ess}\sup\{|F(e^{i\theta})|:\; \theta\in [0, 2\pi)\} &\mbox{if } p=\infty.
\end{cases}
$$

For $\theta\in[0,2\pi]$ and $z\in\mathbb{D}$, let
$$P(z,e^{i\theta})=\frac{1}{2\pi}\frac{1-|z|^{2}}{|1-ze^{-i\theta}|^{2}}$$ be the {\it Poisson kernel}.
For a mapping $F\in L^{1}(\mathbb{T})$, the {\it Poisson integral} of $F$
is defined by

$$f(z)=P[F](z)=\int_{0}^{2\pi}P(z,e^{i\theta})F(e^{i\theta})d\theta.$$

It is well-known that if $F$ is absolutely continuous, then it is of bounded variation. This implies that for
almost all $e^{i\theta}\in\mathbb{T}$, the derivative $\dot{F}(e^{i\theta})$ exists, where
$$\dot{F}(e^{i\theta}):=\frac{dF(e^{i\theta})}{d\theta}.$$

In \cite{Zhu}, the author posed the following problem.

\begin{prob}\label{ques-1}
What conditions on the boundary function $F$ ensure that the partial
derivatives of its harmonic extension $f=P[F]$, i.e., $f_{z}$ and $\overline{f_{\overline{z}}}$, are in the space $\mathcal{B}^{p}(\mathbb{D})$ (or
$\mathcal{H}^{p}(\mathbb{D})$), where $p\geq1$?
\end{prob}

In \cite{Zhu}, the author discussed Problem \ref{ques-1} under the condition that $F$ is absolutely continuous.
First, he proved the following, which is one of the two main results in  \cite{Zhu}. On the related discussion, we refer to the recent paper \cite{K-2019}.

\begin{Thm}{\rm (\cite[Theorem 1.2]{Zhu})}\label{Zhu-1}
Suppose that $p\in[1,2)$ and $f=P[F]$ is a harmonic mapping in $\mathbb{D}$ with
$\dot{F}\in L^{p}(\mathbb{T})$, where $F$ is an absolutely continuous function. Then both $f_{z}$ and $\overline{f_{\overline{z}}}$ are in $\mathcal{B}^{p}(\mathbb{D}).$
\end{Thm}

Furthermore, by requiring the mappings $P[F]$ to be harmonic quasiregular, the interval of $p$ is widened
from $[1,2)$ into $[1,\infty)$, as shown in the following result, which is the other main result in \cite{Zhu}.

\begin{Thm}{\rm (\cite[Theorem 1.3]{Zhu})}\label{Zhu-2}
Suppose that $p\in[1,\infty]$ and $f=P[F]$ is a harmonic $K$-quasiregular mapping in $\mathbb{D}$ with
$\dot{F}\in L^{p}(\mathbb{T})$, where $F$ is an  absolutely continuous function and $K\geq1$. Then
both $f_{z}$ and $\overline{f_{\overline{z}}}$ are in $\mathcal{H}^{p}(\mathbb{D})$.
\end{Thm}

The purpose of this paper is to discuss these two results further. Regarding Theorem \Ref{Zhu-1}, our result is as follows, which shows that Theorem \Ref{Zhu-1} is true for $p\in [1,\infty)$, and also indicates that Theorem \Ref{Zhu-1} is not true when $p=\infty$.

\begin{thm}\label{cp-1.0}
Suppose that $f=P[F]$ is a harmonic mapping in $\mathbb{D}$ and
$\dot{F}\in L^{p}(\mathbb{T})$, where $F$ is an  absolutely continuous  function. \ben
\item\label{cp-1.0-1}
If $p\in[1,\infty)$, then both $f_{z}$ and $\overline{f_{\overline{z}}}$ are in $\mathcal{B}^{p}(\mathbb{D}).$
\item\label{cp-1.0-2}
If $p=\infty$, then there exists a harmonic
mapping $f=P[F]$, where $F$ is an absolutely continuous  function with $\dot{F}\in L^{\infty}(\mathbb{T})$, such that
neither $f_{z}$ nor $\overline{f_{\overline{z}}}$ is in $\mathcal{B}^{\infty}(\mathbb{D}).$
\een
\end{thm}

About Theorem \Ref{Zhu-2}, we show that this result also holds true for harmonic elliptic mappings, which are more general than harmonic quasiregular mappings. In order to state our result, we need to introduce the definition of elliptic mappings.

A sense-preserving continuously differentiable mapping
$f:~\mathbb{D}\rightarrow\mathbb{C}$ is said to be a {\it $(K,K')$-elliptic mapping} if
$f$ is absolutely continuous on lines in $\mathbb{D}$, and there are constants $K\geq1$
and $K'\geq0$ such that
$$\|D_{f}(z)\|^{2}\leq KJ_{f}(z)+K'$$ in $\mathbb{D}.$
In particular, if $K'\equiv0$, then a $(K,K')$-elliptic mapping is said to be {\it $K$-quasiregular}.
It is well known that every quasiregular mapping is an elliptic mapping. But the inverse of this statement is not true. This can be seen from the example:
Let $f(z)=z+\overline{z}^{2}/2$ in $\mathbb{D}$ which is indeed a univalent harmonic mapping of $\mathbb{D}$.
Then elementary computations show that $(a)$ $\ds \sup_{z\in\mathbb{D}}\{|\omega(z)|\}=1$, which implies that $f$ is not
$K$-quasiregular for any $K\geq 1$, and $(b)$ $f$ is a $(1,4)$-elliptic mapping.
We refer  to \cite{C-K,CP-2020,FS,K8,Ni} for more details of elliptic mappings.

Now, we are ready to state our next result.
\begin{thm}\label{cp-1.1}
Suppose that $p\in[1,\infty]$ and $f=P[F]$ is a $(K,K')$-elliptic mapping in $\mathbb{D}$ with
$\dot{F}\in L^{p}(\mathbb{T})$, where $F$ is an absolutely continuous function, $K\geq1$ and $K'\geq0$. Then both
$f_{z}$ and $\overline{f_{\overline{z}}}$ are in $\mathcal{H}^{p}(\mathbb{D}).$
\end{thm}

The proofs of Theorems \ref{cp-1.0} and \ref{cp-1.1} will be presented in Section \ref{csw-sec1.0}.

\section{Proofs of the main results}\label{csw-sec1.0}
We start this section by recalling the following two lemmas from \cite{Zhu}.

\begin{Lem}{\rm (\cite[Theorem 1.1]{Zhu})}\label{Zhu-0.1}
Suppose $p\in[1,\infty)$ and $f=P[F]$ is a harmonic mapping in $\mathbb{D}$ with
$\dot{F}\in L^{p}(\mathbb{T})$, where $F$ is an  absolutely continuous function. Then
for $z=re^{it}\in\mathbb{D}$,
$$\|f_{r}\|_{L^{p}}\leq\big(2C(p)\big)^{\frac{1}{p}}\|\dot{F}\|_{L^{p}},$$
and thus, $f_{r}\in\mathcal{B}^{p}_{\mathcal{G}}(\mathbb{D}),$
where $$C(p)=\int_{0}^{1}\left(\frac{4\tanh^{-1}r}{\pi r}\right)^{p}rdr\leq\frac{4^{p-1}}{\pi^{p}}\big(2^{p}+(2-2^{-p})\Gamma(1+p)\big)$$
and $\Gamma$ denotes the usual Gamma function.
\end{Lem}

\begin{Lem}{\rm (\cite[Lemma 2.3]{Zhu})}\label{Zhu-0.2}
Assume the hypotheses of Lemma \Ref{Zhu-0.1}.
Then for $z=re^{it}\in\mathbb{D}$,
$$\|f_{t}\|_{p}\leq\|\dot{F}\|_{L^{p}},$$
and thus, $f_{t}\in\mathcal{H}^{p}_{\mathcal{G}}(\mathbb{D})$.
\end{Lem}

\subsection{Proof of Theorem \ref{cp-1.0}}
For the proof of the first statement of the theorem, let $z=re^{it}\in\mathbb{D}$. Then we have
\be\label{eq-j-0}f_{t}(z):=\frac{\partial f(z)}{\partial t}=i\big(zf_{z}(z)-\overline{z}f_{\overline{z}}(z)\big)~\mbox{and}~f_{r}(z):=\frac{\partial f(z)}{\partial r}=f_{z}(z)e^{it}+f_{\overline{z}}(z)e^{-it},\ee
which implies that
\beqq\label{eq-j-1}f_{z}(z)=\frac{e^{-it}}{2}\left(f_{r}(z)-\frac{i}{r}f_{t}(z)\right)~\mbox{and}~
\overline{f_{\overline{z}}(z)}=\frac{e^{-it}}{2}\left(\overline{f_{r}(z)}-\frac{i}{r}\overline{f_{t}(z)}\right).\eeqq
It follows that for $p\in[1,\infty)$,
\beqq|f_{z}(z)|^{p}\leq\frac{1}{2^{p}}\left(|f_{r}(z)|+\left|\frac{f_{t}(z)}{r}\right|\right)^{p}
\leq\frac{1}{2}\left(|f_{r}(z)|^{p}+\left|\frac{f_{t}(z)}{r}\right|^{p}\right)
\eeqq
and similarly,
\beqq |f_{\overline{z}}(z)|^{p}\leq \frac{1}{2}\left(|f_{r}(z)|^{p}+\left|\frac{f_{t}(z)}{r}\right|^{p}\right).
\eeqq

Obviously, to prove that $f_{z}$ and $\overline{f_{\overline{z}}}$ are in $\mathcal{B}^{p}(\mathbb{D})$, it suffices to show the following:
\beqq\label{eq-j-2} \int_{\mathbb{D}}|f_{r}(z)|^{p}d\sigma(z)<\infty\;\;\mbox{and}\;\; \int_{\mathbb{D}}\left|\frac{f_{t}(z)}{r}\right|^{p}d\sigma(z)<\infty.\eeqq

We only need to check the boundedness of the integral $\int_{\mathbb{D}}\left|\frac{f_{t}(z)}{r}\right|^{p}d\sigma(z)$ because
the boundedness of the integral $\int_{\mathbb{D}}|f_{r}(z)|^{p}d\sigma(z)$ easily follows from Lemma \Ref{Zhu-0.1}.

By Lemma \Ref{Zhu-0.2}, we have
\beqq \frac{1}{2\pi}\int_{0}^{2\pi}|f_{t}(re^{it})|^{p}dt\leq\|\dot{F}\|_{L^{p}}^{p},\eeqq
which yields that

\be\label{eq-j-4} \int_{\mathbb{D}\backslash\mathbb{D}_{\frac{1}{2}}}\left|\frac{f_{t}(z)}{r}\right|^{p}d\sigma(z)\leq
\frac{2^{p-1}}{\pi}\int_{\frac{1}{2}}^{1}\left(\int_{0}^{2\pi}\left|f_{t}(re^{it})\right|^{p}dt\right)dr\leq2^{p-1}\|\dot{F}\|_{L^{p}}^{p}.\ee

To demonstrate the boundedness of the integral $\int_{\mathbb{D}_{\frac{1}{2}}}\left|\frac{f_{t}(z)}{r}\right|^{p}d\sigma(z)$, assume that $f=h+\overline{g}$ where both $h$ and $g$ being analytic in $\mathbb{D}$. Then $\|D_{f}\|=|h'|+|g'|$. This implies that $\|D_{f}\|$ is continuous in $\overline{\mathbb{D}}_{\frac{1}{2}}$, and thus, $\|D_{f}\|$ is bounded in $\overline{\mathbb{D}}_{\frac{1}{2}}$.
Hence, by (\ref{eq-j-0}), we have
\beq\label{eq-j-5} \int_{\mathbb{D}_{\frac{1}{2}}}\left|\frac{f_{t}(z)}{r}\right|^{p}d\sigma(z)&=&
\int_{0}^{\frac{1}{2}}\int_{0}^{2\pi}r\big|e^{it}f_{z}(re^{it})-e^{-it}f_{\overline{z}}(re^{it})\big|^{p}dt dr\\ \nonumber
&\leq&\int_{0}^{\frac{1}{2}}\int_{0}^{2\pi}r\|D_{f}(re^{it})\|^{p}dt dr\\ \nonumber
&=&\int_{\mathbb{D}_{\frac{1}{2}}}\|D_{f}(z)\|^{p}d\sigma(z)<\infty.
\eeq
Combining (\ref{eq-j-4}) and (\ref{eq-j-5}) gives the final estimate

\beqq
\int_{\mathbb{D}}\left|\frac{f_{t}(z)}{r}\right|^{p}dA(z)=\int_{\mathbb{D}_{\frac{1}{2}}}\left|\frac{f_{t}(z)}{r}\right|^{p}d\sigma(z)
+\int_{\mathbb{D}\backslash\mathbb{D}_{\frac{1}{2}}}\left|\frac{f_{t}(z)}{r}\right|^{p}d\sigma(z)<\infty,
\eeqq
which is what we need, and so, the statement \eqref{cp-1.0-1} of the theorem is true.\medskip

To prove the second statement of the theorem, let $F(e^{i\theta})=|\sin \theta|$, where $\theta\in[0,2\pi]$. Then $F$ is absolutely continuous and
 $\dot{F}\in L^{\infty}(\mathbb{T})$.
Also, elementary computations guarantee that for $z=re^{it}\in\mathbb{D}$,
\beqq f(z)&=&P[F](z)=\int_{0}^{2\pi}P(z,e^{i\theta})|\sin \theta|d\theta\\
&=&\frac{1}{2\pi r(r^{2}-1)}\bigg [(1-r^{2})
\cos t\log\frac{1+r^{2}-2r\cos t}{1+r^{2}+2r\cos t}\\
&&+2(1+r^{2})\sin t
\Big(\arctan\Big(\frac{1+r}{r-1}\cot\frac{t}{2}\Big)+
\arctan\Big(\frac{1+r}{r-1}\tan\frac{t}{2}\Big)\Big)\bigg].
\eeqq
Then $$|f_{z}(z)|=\frac{1}{2}\left|f_{r}(z)-i\frac{f_{t}(z)}{r}\right|=
\frac{1}{2}\sqrt{|f_{r}(z)|^{2}+\frac{|f_{t}(z)|^{2}}{r^{2}}},$$
which implies that
\be\label{twd-1}|f_{z}(r)|=\frac{1}{2}\sqrt{|f_{r}(r)|^{2}+\frac{|f_{t}(r)|^{2}}{r^{2}}}.\ee
Since $$f_{r}(r)=\frac{1}{\pi r^{2}}\log\left(\frac{1-r}{1+r}\right)+\frac{2}{\pi}\frac{1}{r(1-r^{2})},$$
 we see that
\be\label{twd-2}\lim_{r\rightarrow1^{-}}f_{r}(r)=\infty.\ee

Combining  (\ref{twd-1}) and (\ref{twd-2}) gives
$$\lim_{r\rightarrow1^{-}}|f_{z}(r)|=\infty,$$ which implies that $f_{z}$ is not in $\mathcal{B}^{\infty}(\mathbb{D}).$

By the similar reasoning, we know that $\overline{f_{\overline{z}}}$ is not in $\mathcal{B}^{\infty}(\mathbb{D})$ either, and hence,
the theorem is proved.
\qed

\subsection{Proof of Theorem \ref{cp-1.1}} Assume that $f=P[F]$ is a $(K,K')$-elliptic mapping in $\mathbb{D}$, which means that
for $z\in\mathbb{D},$
\be\label{hpw-1}\|D_{f}(z)\|^{2}\leq K\|D_{f}(z)\|l(D_{f}(z))+K'.\ee

We divide the proof of this theorem into two cases.
\bca\label{cla-1}
Suppose that $p\in[1,\infty).$ \eca

 It follows from \eqref{hpw-1} that
\beqq
\|D_{f}(z)\|^{p}&\leq&\bigg(\frac{Kl(D_{f}(z))+\sqrt{\big(Kl(D_{f}(z))\big)^{2}+4K'}}{2}\bigg)^{p}\\ \nonumber
&\leq&
\left(Kl(D_{f}(z))+\sqrt{K'}\right)^{p}\leq2^{p-1}\left(K^{p}l^{p}(D_{f}(z))+K'^{\frac{p}{2}}\right),
 \eeqq
and thus, we have
\be\label{cpw-3.0} l^{p}(D_{f}(z))\geq\frac{1}{2^{p-1}K^{p}}\|D_{f}(z)\|^{p}-\frac{K'^{\frac{p}{2}}}{K^{p}}.\ee

By (\ref{eq-j-0}), (\ref{cpw-3.0}) and Lemma \Ref{Zhu-0.2}, we know that for $z=re^{it}\in\mathbb{D},$

\beqq
2\pi\|\dot{F}\|_{L^{p}}^{p}&\geq&\int_{0}^{2\pi}|f_{t}(re^{it})|^{p}dt\geq r^{p}\int_{0}^{2\pi}l^{p}(D_{f}(re^{it}))dt\\
&\geq&\frac{r^{p}}{2^{p-1}K^{p}}\int_{0}^{2\pi}\|D_{f}(re^{it})\|^{p}dt-\frac{2\pi K'^{\frac{p}{2}}}{K^{p}},
\eeqq
which implies that

\beqq \sup_{r\in(0,1)}\left(\frac{1}{2\pi}\int_{0}^{2\pi}\|D_{f}(re^{it})\|^{p}dt\right)^{\frac{1}{p}}\leq
2^{\frac{p-1}{p}}\left(K^{p}\|\dot{F}\|_{L^{p}}^{p}+K'^{\frac{p}{2}}\right)^{\frac{1}{p}}. \eeqq
Hence $f_{z},~\overline{f_{\overline{z}}}\in \mathcal{H}^{p}(\mathbb{D}).$

\bca\label{cla-2}
Suppose that $p=\infty.$ \eca

By (\ref{hpw-1}), we have

\beqq
\|D_{f}(z)\|\leq\frac{Kl(D_{f}(z))+\sqrt{\big(Kl(D_{f}(z))\big)^{2}+4K'}}{2}
\leq Kl(D_{f}(z))+\sqrt{K'},
 \eeqq
which, together with (\ref{eq-j-0}) and Lemma \Ref{Zhu-0.2}, gives

$$\|\dot{F}\|_{\infty}\geq\|f_{t}\|_{\infty}\geq|f_{t}(re^{it})|\geq rl(D_{f}(re^{it}))\geq\frac{r}{K}\left(\|D_{f}(re^{it})\|-\sqrt{K'}\right).$$
Consequently, $$\sup_{z\in\mathbb{D}}\big(|z|\|D_{f}(z)\|\big)\leq\sqrt{K'}+K\|\dot{F}\|_{\infty},$$
from which we conclude that $f_{z},~\overline{f_{\overline{z}}}\in \mathcal{H}^{\infty}(\mathbb{D}),$ and hence the theorem is proved.
\qed

\normalsize

\end{document}